\documentclass[12pt]{article}
\usepackage{latexsym}
\usepackage{amsfonts}
\usepackage{amsmath}
\usepackage{amssymb}
\usepackage{graphicx}
\usepackage{latexsym}
\usepackage{amsfonts}
\usepackage{graphicx}
\usepackage{psfrag}

\input{tcilatex}
\begin{document}

\qquad %\setcounter{section}{-1}

\thispagestyle{empty}

\begin{center}
{\Large \textbf{\ Means Moments and Newton's Inequalities }}

\vskip0.2inR. Sharma, A. Sharma, R. Saini and G. Kapoor

Department of Mathematics \& Statistics\\[0pt]
Himachal Pradesh University\\[0pt]
Shimla - 171005,\\[0pt]
India \\[0pt]
email: rajesh.sharma.hpn@nic.in
\end{center}

\vskip1.5in \noindent \textbf{Abstract. }It is shown that Newton's
inequalities and the related Maclaurin's inequalities provide several
refinements of the fundamental Arithmetic mean - Geometric mean - Harmonic
mean inequality in terms of the means and variance of positive real numbers.
We also obtain some inequalities involving third and fourth central moments
of real numbers.



\bigskip

\bigskip

\bigskip

\bigskip

\bigskip

\bigskip

\section{\protect\bigskip Introduction}

\setcounter{equation}{0} \ Let $x_{1},x_{2},...,x_{n}$ denote $n$ real
numbers. Their rth moment and rth central moment are respectively the numbers%
\begin{equation}
m_{r}^{\prime }=\frac{1}{n}\overset{n}{\underset{i=1}{\sum }}x_{i}^{r} 
\tag{1.1}
\end{equation}%
and%
\begin{equation}
m_{r}=\frac{1}{n}\overset{n}{\underset{i=1}{\sum }}\left(
x_{i}-m_{1}^{\prime }\right) ^{r},  \tag{1.2}
\end{equation}%
where $r=1,2,\ldots $ $\ \ $. Note that $m_{1}=0,$ $m_{1}^{\prime }$ is the
arithmetic mean and $m_{2}$ is the variance of $n$ real numbers $%
x_{1},x_{2},...,x_{n}.$ It is customary to denote $m_{1}^{\prime }$ and $%
m_{2}$ respectively by $A$ and $s^{2}$. The Geometric mean $\left( G\right) $
and Harmonic mean $\left( H\right) $ of $n$ positive real numbers $%
x_{1},x_{2},...,x_{n}$ are respectively the numbers

\begin{equation*}
G=\left( \overset{n}{\underset{i=1}{\dprod }}x_{i}\right) ^{\frac{1}{n}}%
\text{ and \ }H=\left( \frac{1}{n}\underset{i=1}{\overset{n}{\sum }}\frac{1}{%
x_{i}}\right) ^{-1}.
\end{equation*}%
The $k^{\text{th}}$ elementary symmetric function $C_{k\text{ }}$ of $%
x_{1},x_{2},...,x_{n}$ is the \ sum of the products taken $k$ at a time of
different $x_{i}$'s$.$ The $k^{\text{th}}$ elementary symmetric mean is the
number

\begin{equation}
S_{k}=\frac{1}{\binom{n}{k}}C_{k}.  \tag{1.3}
\end{equation}%
We write $S_{0}=C_{0}=1.$

A well known theorem of Newton (1707) says that elementary symmetric means
satisfy the inequality

\begin{equation}
S_{k}^{2}\geq S_{k+1}S_{k-1}\text{ ,}  \tag{1.4}
\end{equation}%
for all $k=1,2,\ldots ,n,$ \ and equality holds if and only if the numbers $%
x_{i}$'s are all equal.

A related \ theorem due to Maclaurin (1729) asserts that the elementary
symmetric means of $n$ positive real numbers satisfy the following
inequalities string:

\begin{equation}
S_{1}\geq S_{2}^{\frac{1}{2}}\geq \ldots \geq S_{n}^{\frac{1}{n}}\text{ }. 
\tag{1.5}
\end{equation}%
Equality occurs if and only if $x_{i}$'s are all equal. One can easily see
that the inequalities (1.5) also follow from the inequalities (1.4). For
more detail, see Hardy et al (1952).

Newton's identities give the relations between the elementary symmetric
functions and the sum of the powers of $x_{1},x_{2},...,x_{n}.$ Let

\begin{equation*}
\alpha _{k}=\overset{n}{\underset{i=1}{\sum }}x_{i}^{k}\text{ ,}
\end{equation*}%
and write $\alpha _{0}=1.$ Then Newton's identities say that

\begin{equation}
kC_{k}=\underset{i=1}{\overset{k}{\sum }}\left( -1\right)
^{i-1}C_{k-i}\alpha _{i}.  \tag{1.6}
\end{equation}%
The well known Arithmetic mean - Geometric mean - Harmonic mean inequality
(A-G-H inequality) says that $H\leq G\leq A$ . Alternative proofs,
refinements and generalizations of A-G-H inequality have been studied
extensively in literature over roughly the last two centuries. See Bullen
(2003). Note that $S_{1}=A$ and $S_{n}^{\frac{1}{n}}=G.$ The Maclaurin
inequalities (1.5) therefore provide the refinements of the A-G inequality.
These refinements are in terms of the symmetric means $S_{k}.$ On using
Newton's identity we can find the relation between symmetric means and
moments of the $n$ real numbers. Then, the Maclaurin inequalities can be
used to find the\ refinements of the A-G inequality in terms of the
expressions involving means and moments. The purpose of the present note is
to point out some such inequalities which remains unnoticed in the
literature in the present explicit forms.

Our first theorem gives a refinement of the A-G-H inequality in terms \ of
expressions involving arithmetic mean and harmonic mean (see Theorem 2.1,
below). A refinement of the A-G inequality involving second moment is given.
This also provides an upper bound for the variance, and Brunk inequalities
(Theorem 2.2). A further refinement of this inequality involves harmonic
mean (Theorem 2.3). We derive bounds for the third and fourth moments and
central moments of statistical interest (Theorem 2.4 and 2.7). The
refinements of A-G inequality in terms of first four moments are also given
(Theorem 2.5, 2.6 and 2.8).

\section{Main results}

\setcounter{equation}{0}\textbf{Theorem 2.1. }For $n$ positive real numbers $%
x_{1},x_{2},...,x_{n}$ and with notations as above, we have

\begin{equation}
H\leq H\left( \frac{A}{H}\right) ^{\frac{1}{n}}\leq G\leq \left( \frac{H}{A}%
\right) ^{\frac{1}{n}}A\leq A.  \tag{2.1}
\end{equation}%
\textbf{Proof. }The third inequality (2.1) follows from the inequality $%
S_{1}\geq S_{n-1}^{\frac{1}{n-1}}$ in (1.5) and the fact that from (1.3) we
have $S_{1}=A$ and $S_{n-1}H=G^{n}.$ The second inequality (2.1) follows on
applying third inequality (2.1) to the numbers $\frac{1}{x_{1}},\frac{1}{%
x_{2}},\ldots ,\frac{1}{x_{n}}.$ The extreme inequalities are immediate.\ $%
\blacksquare $

It may be noted here that from (2.1) we also have

\begin{equation*}
\left( \frac{G}{H}\right) ^{\frac{1}{n-1}}G\leq A\leq \left( \frac{G}{H}%
\right) ^{n-1}G.
\end{equation*}

and

\begin{equation*}
\left( \frac{G}{A}\right) ^{n-1}G\leq H\leq \left( \frac{G}{A}\right) ^{%
\frac{1}{n-1}}G.
\end{equation*}

\textbf{Theorem 2.2. }For $n$ positive real numbers $x_{1},x_{2},...,x_{n}$
and with notations as above, we have

\begin{equation}
G\leq \sqrt{\frac{nA^{2}-m_{2}^{\prime }}{n-1}}\leq A.  \tag{2.2}
\end{equation}

\textbf{Proof. }It follows from the inequality $S_{2}\geq \left(
S_{n}\right) ^{\frac{2}{n}}$ in (1.5) that for $n$ positive real numbers $%
x_{1},x_{2},...,x_{n}$ the inequality

\begin{equation}
\frac{2}{n\left( n-1\right) }\overset{n}{\underset{i<j}{\sum }}%
x_{i}x_{j}\geq G^{2},  \tag{2.3}
\end{equation}%
holds true. From (1.6), $2C_{2}=C_{1}\alpha _{1}-C_{0}\alpha _{2}$ where $%
C_{0}=1,$ $C_{1}=\alpha _{1}=\sum x_{i}$ $,$ $\alpha _{2}=\sum x_{i}^{2}$
and $C_{2}=\underset{i<j}{\sum }x_{i}x_{j}$\textbf{\ }$.$ Therefore, on
using (1.1), we get that

\begin{equation}
2\underset{i<j}{\overset{n}{\sum }}x_{i}x_{j}=\left( \underset{i=1}{\overset{%
n}{\sum }}x_{i}\right) ^{2}-\overset{n}{\underset{i=1}{\sum }}%
x_{i}^{2}=n\left( nA^{2}-m_{2}^{\prime }\right) .  \tag{2.4}
\end{equation}%
Insert (2.4) in (2.3); a little computation leads to first inequality (2.2).
The second inequality (2.2) follows from the fact that \ $m_{2}^{\prime
}\geq A^{2}.\blacksquare $

The inequality (2.1) can be written in several different ways, for example
the inequality%
\begin{equation*}
G\leq A\sqrt{1-\frac{1}{n-1}\left( \frac{s}{A}\right) ^{2}}
\end{equation*}%
follows from (2.1) on using $m_{2}^{\prime }=s^{2}+A^{2}.$ We also have

\begin{equation}
m_{2}^{\prime }\leq nA^{2}-\left( n-1\right) G^{2}  \tag{2.5}
\end{equation}%
and 
\begin{equation}
s^{2}\leq \left( n-1\right) \left( A^{2}-G^{2}\right) .  \tag{2.6}
\end{equation}%
We note that Brunk's inequalities (1959) follow from the inequality (2.6).
Let $a\leq x_{i}\leq b$ for all $i=1,2,\ldots ,n.$ From (2.6), $s^{2}\leq
\left( n-1\right) A^{2}$. Apply this to $n$ positive numbers $x_{i}-a$, we
get Brunk's inequality, $s\leq \sqrt{n-1}\left( A-a\right) $. Likewise, the
inequality $s\leq \sqrt{n-1}\left( b-A\right) $ follows from (2.6). For the
related variance upper bounds, see Bhatia and Davis (2000), Sharma (2008)
and Sharma et al (2010).

\textbf{Theorem 2.3. }With conditions as in above theorem and for $n\geq 3$,
\ we have

\begin{equation}
G\leq \left( \frac{G^{n}}{H}\right) ^{\frac{1}{n-1}}\leq \sqrt{\frac{%
nA^{2}-m_{2}^{\prime }}{n-1}}\leq A.  \tag{2.7}
\end{equation}%
\textbf{Proof. }From (1.5) , $S_{2}\geq \left( S_{n-1}\right) ^{\frac{2}{n-1}%
}.$ Therefore, for $n$ positive real numbers $x_{1},x_{2},...,x_{n}$ the
inequality

\begin{equation}
\frac{2}{n\left( n-1\right) }\overset{n}{\underset{i<j}{\sum }}%
x_{i}x_{j}\geq \left( \frac{C_{n-1}}{n}\right) ^{\frac{2}{n-1}},  \tag{2.8}
\end{equation}%
holds true. Also

\begin{equation}
C_{n-1}=\overset{n}{\underset{i=1}{\Pi }}x_{i}\underset{i=1}{\overset{n}{%
\sum }}\frac{1}{x_{i}}=n\frac{G^{n}}{H}.  \tag{2.9}
\end{equation}%
Insert (2.4) and \ (2.9) in (2.8); a little computation leads to the second
inequality (2.7) \ $\blacksquare $

From (2.7) , we also have%
\begin{equation*}
G\leq \left( \frac{G^{n}}{H}\right) ^{\frac{1}{n-1}}\leq A\sqrt{1-\frac{1}{%
n-1}\left( \frac{s}{A}\right) ^{2}}
\end{equation*}%
\begin{equation}
m_{2}^{\prime }\leq nA^{2}-\left( n-1\right) \left( \frac{G^{n}}{H}\right) ^{%
\frac{2}{n-1}}.  \tag{2.10}
\end{equation}%
and 
\begin{equation}
s^{2}\leq \left( n-1\right) \left( A^{2}-\left( \frac{G^{n}}{H}\right) ^{%
\frac{2}{n-1}}\right) .  \tag{2.11}
\end{equation}%
The inequalities (2.10) and (2.11) provide refinements of the inequalities
(2.5) and (2.6) , respectively.

\textbf{Theorem 2.4}. For $n$ real numbers $x_{1},x_{2},...,x_{n}$ and $%
n\geq 3$, we have

\begin{equation}
Am_{3}^{\prime }\leq \frac{1}{2\left( n-1\right) }\left[ \left( n-2\right)
m_{2}^{\prime 2}+n\left( n+1\right) m_{2}^{\prime }A^{2}\right] .  \tag{2.12}
\end{equation}%
\textbf{Proof.} For $k=2,$ the inequality (1.4) gives $S_{2}^{2}\geq
S_{1}S_{3}.$ Therefore, for $n$ real numbers $x_{1},x_{2},...,x_{n},$ we have

\begin{equation}
\left( \frac{2}{n\left( n-1\right) }\underset{i<j}{\overset{n}{\sum }}%
x_{i}x_{j}\right) ^{2}\geq \frac{6A}{n\left( n-1\right) \left( n-2\right) }%
\underset{i<j<k}{\overset{n}{\sum }}x_{i}x_{j}x_{k}\text{ }.  \tag{2.13}
\end{equation}%
Further, from (1.6), $3C_{3}=C_{2}\alpha _{1}-C_{1}\alpha _{2}+\alpha _{3}$
where $\alpha _{3}=\sum x_{i}^{3}$ and $C_{3}=\underset{i<j<l}{\sum }%
x_{i}x_{j}x_{l}$ $.$ Therefore, on using (2.3) and (1.1), we find that

\begin{equation*}
3\underset{i<j<k}{\overset{n}{\sum }}x_{i}x_{j}x_{k}=\underset{i=1}{\overset{%
n}{\sum }}x_{i}^{3}+\frac{1}{2}\left( \underset{i=1}{\overset{n}{\sum }}%
x_{i}\right) ^{3}-\frac{3}{2}\underset{i=1}{\overset{n}{\sum }}x_{i}\underset%
{i=1}{\overset{n}{\sum }}x_{i}^{2}
\end{equation*}

\begin{equation}
=n\left( m_{3}^{\prime }+\frac{n^{2}}{2}A^{3}-\frac{3}{2}nAm_{2}^{\prime
}\right) .  \tag{2.14}
\end{equation}%
On inserting (2.4) and (2.14) in (2.13) and simplifying the resulting
expressions; we immediately get (2.12). \ $\blacksquare $

The twin inequality

\begin{equation}
Am_{3}\leq \frac{n-2}{2\left( n-1\right) }s^{4}+\frac{\left( n-2\right) }{2}%
A^{2}s^{2}  \tag{2.15}
\end{equation}%
follows from (2.12) on using the relations $m_{2}^{\prime }=s^{2}+A^{2}$ and 
$m_{3}^{\prime }=m_{3}+3As^{2}+A^{3}$ $.$It may be noted here that

\begin{equation*}
3\underset{i<j<k}{\overset{n}{\sum }}x_{i}x_{j}x_{k}=n\left( m_{3}-\frac{%
3\left( n-2\right) }{2}As^{2}+\frac{\left( n-1\right) \left( n-2\right) }{2}%
A^{3}\right)
\end{equation*}%
and the inequality

\begin{equation*}
m_{3}\leq \frac{n-2}{2\left( n-1\right) }\frac{s^{4}}{A}+\frac{\left(
n-2\right) }{2}As^{2}
\end{equation*}%
holds true for positive real numbers $x_{1},x_{2},...,x_{n}.$

\ \textbf{Corollary 2.1. }Under the conditions of Theorem 2.2, we have

\begin{equation}
m_{3}\leq \left( n-2\right) As^{2}\leq \left( n-1\right) \left( n-2\right)
A^{3}\text{.}  \tag{2.16}
\end{equation}%
\textbf{Proof. }The inequality (2.16) follows from (2.15) and the fact that
for positive real numbers $x_{1},x_{2},...,x_{n},$ we have $A>0$ and $s\leq 
\sqrt{n-1}$ $A.$ \ $\blacksquare $

Note that the second inequality (2.16) implies that for $n$ real numbers $%
a\leq x_{i}\leq b$, $i=1,2,\ldots ,n$, \ we have $m_{3}\leq \left(
n-1\right) \left( n-2\right) \left( b-A\right) ^{3}$ and $m_{3}\leq \left(
n-1\right) \left( n-2\right) \left( A-a\right) ^{3}.$ Thus, we get
inequalities analogous to Brunk's inequalities (1959) for third central
moment. Likewise, from the first inequality (2.16), we have $m_{3}\leq
\left( n-2\right) \left( b-A\right) s^{2}$ and $m_{3}\leq \left( n-2\right)
\left( A-a\right) s^{2}$. Also, see Sharma et al (2012).

\textbf{Theorem 2.5. }With conditions as in Theorem 2.3 , \ we have%
\begin{equation}
G\leq A\left( \frac{n^{2}}{\left( n-1\right) \left( n-2\right) }\left( 1-%
\frac{3}{n}\frac{m_{2}^{\prime }}{A^{2}}+\frac{2}{n^{2}}\frac{m_{3}^{\prime }%
}{A^{3}}\right) \right) ^{\frac{1}{3}}.  \tag{2.17}
\end{equation}

\textbf{Proof.} It follows from the inequality $S_{3}\geq S_{n}^{\frac{3}{n}%
} $ in (1.5) that for $n$ positive real numbers $x_{1},x_{2},...,x_{n},$ we
have%
\begin{equation}
\frac{6}{n\left( n-1\right) \left( n-2\right) }\underset{i<j<k}{\overset{n}{%
\sum }}x_{i}x_{j}x_{k}\geq G^{3}.  \tag{2.18}
\end{equation}

Insert (2.14) in (2.18); a little computation leads to (2.17). \ $%
\blacksquare $

From (2.17), we also have

\begin{equation*}
G\leq A\left( 1-\frac{3}{n-1}\left( \frac{s}{A}\right) ^{2}+\frac{2}{\left(
n-1\right) \left( n-2\right) }\frac{m_{3}}{A^{3}}\right) ^{\frac{1}{3}}
\end{equation*}%
\begin{equation*}
m_{3}^{\prime }\geq \frac{1}{2}\left[ 3nAm_{2}^{\prime }-n^{2}A^{3}+\left(
n-1\right) \left( n-2\right) G^{3}\right]
\end{equation*}%
and

\begin{equation*}
m_{3}\geq \frac{\left( n-1\right) \left( n-2\right) }{2}\left[ \frac{3}{%
\left( n-1\right) }As^{2}-\left( A^{3}-G^{3}\right) \right] .
\end{equation*}

\textbf{Theorem 2.6.}With conditions as in above theorem and for $n\geq 4$,
\ we have

\begin{equation}
G\leq \left( \frac{G^{n}}{H}\right) ^{\frac{1}{n-1}}\leq A\left( \frac{n^{2}%
}{\left( n-1\right) \left( n-2\right) }\left( 1-\frac{3}{n}\frac{%
m_{2}^{\prime }}{A^{2}}+\frac{2}{n^{2}}\frac{m_{3}^{\prime }}{A^{3}}\right)
\right) ^{\frac{1}{3}}.  \tag{2.19}
\end{equation}

\textbf{Proof.} It follows from inequality $S_{3}\geq S_{n-1}^{\frac{3}{n-1}%
} $ in (1.5) that for $n$ positive real numbers $x_{1},x_{2},...,x_{n},$ we
have%
\begin{equation}
\frac{6}{n\left( n-1\right) \left( n-2\right) }\underset{i<j<k}{\overset{n}{%
\sum }}x_{i}x_{j}x_{k}\geq \left( \frac{G^{n}}{H}\right) ^{\frac{3}{n-1}}. 
\tag{2.20}
\end{equation}

Insert (2.14) in (2.20); a little computation leads to (2.19). \ $%
\blacksquare $

From (2.19), we also have%
\begin{equation*}
G\leq \left( \frac{G^{n}}{H}\right) ^{\frac{1}{n-1}}\leq A\left( 1-\frac{3}{%
n-1}\left( \frac{s}{A}\right) ^{2}+\frac{2}{\left( n-1\right) \left(
n-2\right) }\frac{m_{3}}{A^{3}}\right) ^{\frac{1}{3}}\text{,}
\end{equation*}%
\begin{equation*}
m_{3}^{\prime }\geq \frac{1}{2}\left[ 3nAm_{2}^{\prime }-n^{2}A^{3}+\left(
n-1\right) \left( n-2\right) \left( \frac{G^{n}}{H}\right) ^{\frac{3}{n-1}}%
\right]
\end{equation*}%
and%
\begin{equation*}
m_{3}\geq \frac{\left( n-1\right) \left( n-2\right) }{2}\left[ \frac{3}{%
\left( n-1\right) }As^{2}-\left( A^{3}-\left( \frac{G^{n}}{H}\right) ^{\frac{%
3}{n-1}}\right) \right] \text{.}
\end{equation*}%
Recently , Sharma and Bhandari (2015) have shown that the Newton inequality
provides an inequality involving skewness and kurtosis, 
\begin{equation*}
1+\frac{m_{3}^{2}}{m_{2}^{3}}\leq \frac{m_{4}}{m_{2}^{2}}\leq \frac{1}{2}%
\frac{n-3}{n-2}\frac{m_{3}^{2}}{m_{2}^{3}}+\frac{n}{2}\text{, \ \ }n\geq 3%
\text{.}
\end{equation*}%
Sharma et al (2015) have obtained bounds for the fourth central moment. Here
we show that Newton inequality provide several other inequalities involving
first four moments.

\textbf{Theorem 2.7.} For $n$ real numbers $x_{1},x_{2},...,x_{n}$ and $%
n\geq 4$, we have%
\begin{equation}
\left( nA^{2}-m_{2}^{\prime }\right) m_{4}^{\prime }\geq \frac{1}{6\left(
n-3\right) }\left( 
\begin{array}{c}
n^{4}A^{6}-n^{3}\left( n+4\right) m_{2}^{\prime }A^{4}+4n^{2}\left(
n-1\right) m_{3}^{\prime }A^{3}+9n^{2}m_{2}^{\prime 2}A^{2} \\ 
+4n\left( n-5\right) m_{2}^{\prime }m_{3}^{\prime }A-3n\left( n-2\right)
m_{2}^{\prime 3}-4\left( n-3\right) m_{3}^{\prime 2}%
\end{array}%
\right) \text{.}  \tag{2.21}
\end{equation}

\textbf{Proof. }The inequality (2.21) follows from the inequality $%
S_{3}^{2}\geq S_{2}S_{4}$\textbf{\ }in (1.4). Note that from (1.6), $%
4C_{4}=C_{3}\alpha _{1}-C_{2}\alpha _{2}+C_{1}\alpha _{3}-C_{0}\alpha _{4}$.%
\textbf{\ }Therefore, on using (2.3), (2.14) and (1.1), we get that 
\begin{equation*}
24\underset{i<j<k<l}{\overset{n}{\sum }}%
x_{i}x_{j}x_{k}x_{l}=n^{4}A^{4}-6n^{3}m_{2}^{\prime }+8n^{2}Am_{3}^{\prime
}+3n^{2}m_{2}^{\prime 2}-6nm_{4}^{\prime }\text{ }.
\end{equation*}%
\ $\blacksquare $

\textbf{Theorem 2.8}. For $n$ positive real numbers $x_{1},x_{2},...,x_{n}$
\ and $n\geq 4$, we have

\begin{equation}
G^{4}\leq \frac{1}{\left( n-1\right) \left( n-2\right) \left( n-3\right) }%
\left( n^{3}m_{1}^{\prime 4}-6n^{2}m_{2}^{\prime }m_{1}^{\prime
2}+8nm_{3}^{\prime }m_{1}^{\prime }+3nm_{2}^{\prime 2}-6m_{4}^{\prime
}\right)  \tag{2.22}
\end{equation}%
\textbf{Proof. }From (1.5), $S_{4}\geq \left( S_{n}\right) ^{\frac{4}{n}}.$
Therefore, for $n$ positive real numbers $x_{1},x_{2},...,x_{n}$ the
inequality

\begin{equation*}
\frac{24}{n\left( n-1\right) \left( n-2\right) \left( n-3\right) }\overset{n}%
{\underset{i<j<k<l}{\sum }}x_{i}x_{j}x_{k}x_{l}\geq G^{4},
\end{equation*}%
holds true. \bigskip\ \ $\blacksquare $

The inequality (2.22) also gives%
\begin{equation*}
G\leq A\left( 1-\frac{6}{n-1}\left( \frac{s}{A}\right) ^{2}+\frac{8}{\left(
n-1\right) \left( n-2\right) }\frac{m_{3}}{A^{3}}+\frac{3\left(
ns^{4}-2m_{4}\right) }{\left( n-1\right) \left( n-2\right) \left( n-3\right)
A^{4}}\right) ^{\frac{1}{4}}\text{.}
\end{equation*}

\bigskip \vskip0.2in\noindent \textbf{Acknowledgements}. The authors are
grateful to Prof. Rajendra Bhatia for the useful discussions and
suggestions, and I.S.I. Delhi for a visit in February 2017 when this work
was done. The support of the UGC-SAP is also acknowledged.

\end{document}